\documentclass[11pt]{amsart}
\usepackage{amsfonts,amsmath,amssymb}
\usepackage[all]{xy}

\begin{document}

\newcommand{\PGl}{{\rm PGl}}
\newcommand{\SU}{{\rm SU}}
\newcommand{\U}{{\rm U}}
\newcommand{\cH}{{\mathcal H}}
\newcommand{\vol}{{\rm vol}}
\newcommand{\bS}{{\bf S}}
\newcommand{\op}{{\rm op}}
\newcommand{\dR}{{\rm dR}}
\newcommand{\ie}{\textit{i}.\textit{e}.\,}
\newcommand{\eg}{\textit{e}.\textit{g}.\,}
\newcommand{\cf}{{\textit{cf}.\,}}
\newcommand{\cL}{{\mathcal L}}
\newcommand{\End}{{\rm End}}
\newcommand{\Hom}{{\rm Hom}}
\newcommand{\h}{{\mathfrak h}}
\newcommand{\Z}{{\mathbb Z}}
\newcommand{\T}{{\mathbb T}}
\newcommand{\MU}{{\rm MU}}
\newcommand{\C}{{\mathbb C}}
\newcommand{\B}{{\mathcal B}}
\newcommand{\CP}{{\mathbb{C}P}}
\newcommand{\CPi}{{\mathbb{C}P_\infty}}
\newcommand{\pt}{{\rm pt}}
\newcommand{\Q}{{\mathbb Q}}
\newcommand{\R}{{\mathbb R}}

\title{Formal groups and geometric quantization} 

\author[Jack Morava]{Jack Morava}

\address{Department of Mathematics, The Johns Hopkins University,
Baltimore, Maryland 21218}

\email{jack@math.jhu.edu}


\date{January 2020}

\begin{abstract}{The complex projective spaces, considered as prequantized
symplectic manifolds, are roughly to the complete symmetric functions as
those projective spaces, regarded as complex-oriented manifolds, are to 
Newton's power sums.} \end{abstract}

\maketitle 

{\bf Introduction} \bigskip

This paper is concerned with the statistical mechanics of a commutative
algebra $\B_* \subset \B_* \otimes \R$ of cobordism classes of 
geometrically prequantized symplectic manifolds, which is to some extent
well-understood [3,16]: $\B_*$  is naturally isomorphic to the bordism
ring $\MU_*\CPi$ of compact complex-oriented manifolds carrying a 
suitable complex line bundle with connection, and the principal technical 
result below is an explicit formula [Prop 3.1] for the (injective) Hurewicz 
homomorphism
\[
\h : \B_* \cong \MU_*\CPi \to H_*(\CPi;H_*(\MU,\Z)) \cong \bS_* \otimes 
\Z[b_{(n)} \:|\:n \geq 1]
\]
which identifies these manifolds in terms of their characteristic 
numbers\begin{footnote}{Here $\bS_*$ is the ring of classical symmetric 
functions, and $b_{(n)}$ is the $n$th divided power of an element $b$. Further
technical definitions will be provided soon.}\end{footnote}. This formula is 
easily accessible by modern homotopy theory, though it may not be so familiar; 
the point is that in our case this homomorphism has an interesting 
interpretation as a kind of partition function, analogous to the construction 
which assigns to a Riemannian manifold $(M,g)$, the trace of its 
heat kernel $\exp(t\Delta_g)$ [1,2,5,8,12,21 \dots]. An essentially equivalent
corollary is that (as observed by Friedrich and McKay [6,7]), 
Mi\v{s}\v{c}enko's logarithm
\[
\log_\MU(z) = \sum_{k \geq 1} \frac{\CP_{k-1}}{k} z^k \in MU^2(\CPi)
\]
for complex cobordism can be regarded as a kind of cumulant generating 
function (or Helmholtz free energy, or as (the negative of) a kind of Shannon
entropy) for $\B_*$. \bigskip

\newpage

{\bf \S I, algebraic and geometric preliminaries} \bigskip

{\bf 1.1} For our purposes, a geometrically prequantized symplectic manifold 
[4] will be a $2n$-dimensional closed compact smooth manifold $(V,L,\nabla_L,
j)$ together with a Hermitian complex line bundle $L$ with connection 
$\nabla_L$ on $V$, with symplectic (\ie closed and nondegenerate) connection 
form $\omega = j \Omega(\nabla_L)$ (with $[\omega/2\pi i] \in H^2_\dR(V)$ 
equal to the Chern class $c_1(L)$ of the line bundle $L : V \to \CPi$ (or, 
equivalently, classified by a map to the Narasimhan-Ramanan space $B\T$ of 
bundles with Hermitian connection)). There is a contractible space of 
almost-complex structures $\{j \in \End(T_V) \:|\: j^2 = -1 \}$ on $V$, 
compatible with $\omega$ in the sense that $\omega(j-,-)$ is a Riemannian 
metric, and we assume that such a $j$ has been chosen. The associated 
Liouville volume $\omega^n/n!$ defines a class in $H^{2n}_\dR(V)$ Poincar\'e 
dual to the orientation class in $H_{2n}(V;\Z)$. This will usually be 
summarized as the assertion that $(V,\omega)$ is a (prequantized) symplectic 
manifold. Relaxing the integrality condition on $\omega$ defines the real 
completion $(\B \otimes \R)_*$ of $\B_*$.

Following VL Ginzburg [8,9 Th. H.10], a cobordism
\[
W : (V_0,\omega_0) \to (V_1,\omega_1)
\]
between two such manifolds is a compact $(2n+1)$-dimensional manifold $W$ 
together with a Hermitian line bundle $\cL : W \to B\T$ with connection, 
which restricts to the given line bundle with connection on $\partial W = 
V_0^\op \coprod V_1$. We furthermore assume that the kernel of the 
curvature form $\omega_\cL : T_W \to T^*_W$ is a real line bundle on $W$, 
trivial on $\partial W$. This defines the (monoidal) cobordism category of 
geometrically prequantized manifolds, whose cobordism ring is known:

{\bf Theorem}: The map
\[
\B_* \ni [V,L,\nabla_L,j] \to [V,j,L] \in \MU_*\CPi
\]
defines an isomorphism [3,16] of (evenly graded, commutative) torsion-free 
algebras. \bigskip

{\bf Remark}: Data of the sort described above can be used to define various
(almost K\"ahler, ${\rm Spin}^c$ Dirac, twisted signature \dots) elliptic 
differential operators, which behave nicely with respect to products and 
cobordisms; the work [12] of Liu and Xu, for example, was one of the main
motivations for this paper. Symplectomorphism groups are generally 
infinite-dimensional, however, so it seems unreasonable to expect a close
analog of the rich spectral theory of the classical Laplace-Beltrami operator
in this context. But even if we may not be able to hear the shape of a 
prequantized manifold, we can at least hear its cobordism class. 

{\bf 1.2.1} It will useful to have some examples. We will regard the complex
projective spaces $\CP_n = (\C^\times - 0)/\C^\times$ as symplectic manifolds
$\CP_n(\omega)$ when endowed with the Fubini-Study symplectic form 
\[
\omega = \frac{i}{2} \partial \overline{\partial} \log |z|^2 \in \Omega^2
(\CP_n)\;,
\]
with canonical line bundle pulled back along the standard inclusion 
\[
i_n :\CP_n \to \cup_{n \geq 0} \CP_n := \CPi \cong B\T \;.
\]
More generally, $\CP_n(m\omega)$ will be defined by the line bundle 
$L^{\otimes m}$ with symplectic form $m\omega, \; m \in \Z^\times>0$
\begin{footnote}{We will try to be careful with gradings, which can play a 
rather subtle role. If $A_*$ is a $\Z$-graded module, then $A^* = A_{-*}$. 
We will use graded book-keeping indeterminates, \eg $z_0,z_1,\dots$ of 
cohomological degree +2.}\end{footnote}. Let 
\[
b^\MU_n = [i_n : \CP_n \to \CPi] \in \MU_{2n}\CPi
\]
and define the generating function $b^\MU(z) = 1  + \sum_{n \geq 1} b^\MU_n 
z^n \;,\eg b^\MU_0 = 1$. 

{\bf 1.2.2} Since $\CPi \simeq B\T$ is a (homotopy-commutative) Eilenberg-Mac 
Lane space of type $H(\Z,2)$, 
\[
H_*(\CPi;\Z) \cong \Z[b_{(n)} \:|\: n \geq 1] := \Z[b_{(*)}]
\]
is a bicommutative Hopf algebra under Pontrjagin multiplication, Cartier dual
to the primitively generated Hopf algebra $H^*(\CPi;\Z) \cong \Z[c_1]$. This 
duality implies the relation $b(z_0) \cdot b(z_1) = b(z_0 + z_1)$, where 
$b(z) = 1 + \sum_{n \geq 1} b_{(n)}z^n$ has formal degree zero in 
$H_*(\CPi;\Z[z])$; alternatively, $b_{(n)} = b_{(1)}^n/n!$ is a divided power. 

Steenrod's cycle class homomorphism
\[
h : MU_*X \to H_*(X,\Z)
\]
sends a bordism class $\varkappa = [x : M \to X] \in \MU_*(X)$ to $x_*[M]$,
where $[M] \in H_{2n}(M;\Z)$ is the orientation class. Thus
\[
h(b^\MU_n) = i_{n*}[\CP_n] \in H_{2n}(\CP_\infty;\Z) \;.
\]
The composition 
\[
\xymatrix{
H_*(M;\Z) \ar[r] & H_*(M;\R) \ar[r]^-{(\cap M)^{-1}} & H^{2n-*}_\dR(M) 
\ar[r]^-{*_g} & H^*_\dR(M) }
\]
sends $[M] \in H_{2n}(M;\Z)$ to the Riemannian (hence, in our case, Liouville)
volume 
\[
(\frac{\omega^n}{n!})[M] = \vol(M,\omega) \;;
\]
thus the cycle map sends $\varkappa$ to Weyl's leading term in the asymptotic 
expansion [21] for the trace of the heat kernel $\exp(t \Delta_g)$ of
$(M,\omega)$. \bigskip

{\bf \S II Characteristic numbers} \bigskip

{\bf 2.1} Ravenel and Wilson [19] describe $\MU_*\CPi$ as follows:

{\bf Theorem} {\it As an algebra,
\[
\MU_*\CPi \cong \MU_*[b^\MU_n \:|\: n \geq 1]/(b^\MU(z_0) \cdot b^\MU(z_1)
= b^\MU(z_0 +_\MU z_1)) \;.
\]
It is also a cocommutative Hopf $\MU_*$-algebra, with coproduct} $\Delta b^\MU
(z) = (b^\MU \otimes 1)(z) \otimes_\MU (1 \otimes b^\MU)(z)$.

This is closely related to work [11] of Katz. Here
\[
z_0 +_\MU z_1 = \exp_\MU(\log_\MU(z_0) + \log_\MU(z_1)) 
\]
(with $\exp_\MU$ the formal inverse of $\log_\MU$) is Quillen's formal group
law for complex cobordism; recall that $\MU_*(\pt) = \MU_*$ is polynomial 
over $\Z$ with one generator of each even degree [15,18], and that 
(following  work of Thom) it is generated over $\Q$ (but not $\Z$) by the 
classes $[\CP_n]$. $\Box$

{\bf 2.2} A generator $b_{(1)} : S^2 = \CP_1 \to \CPi$ of $\pi_2\CPi$ defines
a homotopy associative map from the free topological monoid generated by the 
two-sphere, to $\CPi$, for example with the Segre product. Work of IM James 
[20] shows this map to be stably equivalent to a map from $\Omega S^3 \to 
\CPi$. In fact a level one projective representation of the loop group 
L$\SU(2)$ on a separable Hilbert space $\cH$ defines a continuous homomorphism 
\[
\Omega \SU(2) \to \PGl_\C(\cH) \; (\simeq \CPi)
\]
[17] inducing a homomorphism 
\[
\MU_*\Omega \SU(2) \cong \MU_*[b] \to \MU_*[b^\MU_n \:|\: n \geq 1] 
\cong \MU_*\CPi 
\]
of Hopf algebras, taking the class of the adjoint $b : S^2 \to \Omega S^3$ 
of the identity map $S^1 \wedge S^2 \to S^3$ to $b_{(1)}$. Similarly, the 
embedding
\[
\MU^*\CPi = \MU^*[[c]] \to \MU^*\Omega \SU(2) \cong \MU^*[[c_{(n)} \:|\: n 
\geq 1]]
\]
(of the formal group on the left into the completed ring of divided powers on
the right) represents an analog of an exponential map for a formal Lie group. 
From now on we will identify $b_{(1)}$ and $b$. 

{\bf 2.3} The Hurewicz homomorphism 
\[
\h : \MU_*X \cong \pi_*(X \wedge \MU) \to H_*(X \wedge \MU;\Z) \cong H_*(X;\Z) 
\otimes \bS_* 
\]
(where $H_*(\MU;\Z) \cong H_*(B\U;\Z) \cong \bS_*$ is regarded as the 
classical algebra $\Z[h_i \:|\: i \geq 1]$ of complete symmetric functions) 
can be identified, when $H_*(X;\Z)$ is torsion-free, with the map which sends 
$\varkappa = [x : M \to X]$ to 
\[
\h(\varkappa) \in \Hom(H^*(B\U),H_*(X)) \cong H_*(X) \otimes H_*(B\U)
\]
defined by
\[
\h(\varkappa)(\alpha) = x_*D_M\nu^*(\alpha) \;;
\]
where $\nu : M \to B\U$ classifies the stable normal bundle of $M$, and 
$D_M: H^*(M;\Z) \to H_{2n-*}(M;\Z)$ is the Poincar\'e duality 
map\begin{footnote}{From here on, integral (co)homology coefficients will 
often be omitted. The Chern classes of the tangent bundle 
(up to a sign) equal the classes defined by complete symmetric functions of 
roots of the normal bundle $\nu$. Chern-Weil theory presents these global
invariants in terms of local curvature forms.}\end{footnote}.

Following Thom and Milnor, the Hurewicz homomorphism $\MU_*(\pt) \to 
H_*(B\U;\Z)$ is injective, and Quillen's work implies that the image of the 
formal group law on $\MU^*\CPi$ is isomorphic to the additive group over 
$H^*(\CPi;H_*(B\U))$; see [14] for an elegant account. It follows in 
particular that the characteristic number class $\h(\CP_{k-1})$ is divisible 
by $k$. 

Composition with the morphism $[1 : \MU \to H\Z] \in H^0(\MU;\Z)$ of spectra 
yields Steenrod's cycle map
\[
\xymatrix{
h : \MU_*\CPi \ar[r]^\h & H_*\CPi \otimes \bS_* \ar[r]^\varepsilon & H_*\CPi}
\]
($\varepsilon(h_i) = 0, \; i > 0$) sending $b^\MU_n \to i_{n*} D_{\CP_n}
(c_1^0) = b_{(n)}$; more generally, we can think of $h$ as sending $(V,
\omega)$ to its Liouville volume. \bigskip

{\bf \S III Conclusion and final remarks}\bigskip

{\bf 3.1} To state the result below we need some notation for partitions $\pi =
1^{r_1}2^{r_2}\dots$ of 
\[
n = |\pi| = \sum_{k \geq 1} kr_k \;.
\]
We write $r_*$ for the vector $(r_1,r_2,\dots)$ of repetitions in $\pi$ and 
$r = r(\pi) = \sum_{k \geq 1}r_k$ for their sum, and 
\[
\binom{r}{r_*} = \frac{r(\pi)!}{\prod r_k!}
\]
for the associated multinomial function. For example, $r_* = n,0,\dots 
\Rightarrow r=n, \; r_* = 0,\dots,0,1 \Rightarrow r = 1$. \bigskip

{\bf Proposition} 
\[
b^\MU(z) = \exp (b \log_\MU(z)) \in ((\MU \otimes \Q)_*\CPi)[[z]]
\]
{\it and hence}
\[
\h(b^\MU_n) = \sum_{|\pi|=n}\binom{r}{r_*} \: \prod_{k \geq 1} (k^{-1} 
\h(\CP_{k-1}))^{r_k} \: b_{(r)} \in \bS_* \otimes \Z[b_{(*)}] \;. 
\]\bigskip

{\bf Proof} Let $\kappa(z) = \log b^\MU(z) = bz + \dots$, where $\log (1 - x) 
= - \sum_{k \geq 1} z^k/k$ : then $\kappa(z_0 +_\MU z_1) = \kappa(z_0) + 
\kappa(z_1)$. If now $z_i = \exp_\MU(w_i) = w_i + \dots, \; i = 0,1$, then 
$\exp_\MU(w_0 + w_1) = z_0 +_\MU z_1$, so 
\[
\kappa \circ \exp_\MU(w_0 + w_1) = \kappa \circ \exp_\MU(w_0) +  \kappa 
\circ \exp_\MU(w_1)
\]
and hence $\kappa \circ \exp_\MU (w) = bw, \ie b^\MU(\exp_\MU(w)) = \exp(bw)$,
so 
\[
b^\MU(z) = \exp(b\log_\MU(z)) = \exp(b\sum_{k \geq 1} \frac{\CP_{k-1}}{k}z^k) \;;\]
but it is elementary [13 I \S 2.14] that this maps to 
\[ 
\sum_{n \geq 1, |\pi|= n} \binom{r}{r_*} \: 
\prod \frac{(k^{-1}\h(\CP_{k-1}))^{r_k}}{r_k!} \cdot (\sum r_k)! \cdot 
\frac{b_{\sum r_k}}{(\sum r_k)!} \cdot z^n \in (\bS_*[b_{(*)}])[[z]] \;.
\] \bigskip

This can be reformulated as a \bigskip

{\bf Corollary} If $t \in \R_+$ then 
\[
\CP_n(t\omega) = t^n \cdot \vol(\CP_n,\omega) + \dots + t \cdot n^{-1} 
\CP_{n-1} \cdot \vol(\CP_1,\omega) \in (\B \otimes \R)_* 
\]
(analogous to a heat kernel expansion). $\Box$ \bigskip

{\bf 3.2} Such relations are familiar from the theory of symmetric functions. 
If
\[
E(z) = \prod_{i \geq 1} (1 + x_i t) = \sum_{k \geq 0} e_k t^k 
\]
\[
H(z) = \prod_{i \geq 1} (1 - x_i t)^{-1} = \sum_{k \geq 0} h_k t^k
\]
then $E(z)H(-z)^{-1} = 1$, while 
\[
H'(t)/H(t) = \sum_{k \geq 1} p_k t^k
\]
with $p_k = \sum_{i \geq 1} x_i^k$, so 
\[
H(t) = \exp(\sum_{k \geq 1} \frac{p_k}{k} \: t^k) \;.
\]
This suggests a formal analogy in which the symplectic cobordism class of
$\CP_k(\omega)$ is to the complete symmetric function $h_k$ as $b \CP_{k-1}$
is to the power sum $p_k$, with $b$ playing the role of the inverse temperature
$\beta = 1/kT$ in statistical mechanics. \bigskip

\bibliographystyle{amsplain}

\end{document}